\documentclass[a4paper, 12pt]{article}

\usepackage{amsfonts, amsmath, amsthm}
\usepackage{amssymb}
\usepackage{latexsym}
\usepackage[dvips]{graphicx}
\usepackage{color}

\theoremstyle{plain}
\newtheorem{thm}{Theorem}[section]

\newtheorem{lem}{Lemma}[section]

\theoremstyle{definition}

\theoremstyle{definition}
\newtheorem{rmk}{Remark}[section]

\numberwithin{equation}{section}

\newcommand{\R}{\mathbb{R}}
\newcommand{\C}{\mathbb{C}}

\newcommand{\cc}[1]{\overline{#1}}
\newcommand{\op}[1]{\mathcal{#1}}

\newcommand{\pa}{\partial}
\newcommand{\eps}{\varepsilon}
\newcommand{\jb}[1]{\langle #1 \rangle}

\DeclareMathOperator{\realpart}{\rm Re}
\DeclareMathOperator{\imagpart}{\rm Im}

\newcommand{\dis}{\displaystyle}

\begin{document}
\title{
On the derivative nonlinear Schr\"odinger equation 
with weakly dissipative structure\\
 }

\author{
          {Chunhua Li} \thanks{
              Department of Mathematics, College of Science,
              Yanbian University.
              977 Gongyuan Road, Yanji,
              Jilin 133002, China.
              (E-mail: {\tt sxlch@ybu.edu.cn})
             }
          \and
            Yoshinori Nishii \thanks{
              Department of Mathematics, Graduate School of Science,
              Osaka University.
              1-1 Machikaneyama-cho, Toyonaka,
              Osaka 560-0043, Japan.
              (E-mail: {\tt y-nishii@cr.math.sci.osaka-u.ac.jp})             }
           \and
          Yuji Sagawa \thanks{
             Micron Memory Japan, G.K.
             7-10 Yoshikawakogyodanchi, Higashihiroshima,
             Hiroshima 739-0153, Japan. 
             }
           \and
          Hideaki Sunagawa \thanks{
              Department of Mathematics, Graduate School of Science, 
                Osaka City University. 
                3-3-138 Sugimoto, Sumiyoshi-ku, 
                Osaka 558-8585, Japan. 
             }
}

\date{\today }
\maketitle


\noindent{\bf Abstract:}\ We consider the initial value problem for 
cubic derivative nonlinear Schr\"odinger equation in one space dimension. 
Under a suitable weakly dissipative condition on the nonlinearity, 
we show that the small data solution has a logarithmic time decay
in $L^2$.
\\

\\



\section{Introduction and the main result}  \label{sec_intro}
We consider the initial value problem 
\begin{align}
 \left\{\begin{array}{cl}
  i\pa_t u +\frac{1}{2}\pa_x^2u = N(u,\pa_x u), & t>0,\ x \in \R,\\
  u(0,x)=\varphi(x), &x \in \R,
\end{array}\right.
\label{nls}
\end{align}
where $i=\sqrt{-1}$, $u=u(t,x)$ is a $\C$-valued unknown function on 
$[0,\infty)\times\R$.  
$\varphi$ is  a prescribed $\C$-valued function on $\R$ 
which belongs to $H^3\cap H^{2,1}$ and is suitably small in its norm. 
Here and later on as well, for non-negative integers $k$ and $m$, 
$H^k$ denotes the standard $L^2$-based Sobolev 
space of order $k$, and the weighted Sobolev space $H^{k,m}$ is defined 
by $\{\phi \in L^2\,|\, \jb{\, \cdot \,}^{m} \phi \in H^k\}$, 
equipped with the norm 
$\|\phi\|_{H^{k,m}}=\|\jb{\, \cdot \,}^{m} \phi \|_{H^k}$, 
where $\jb{x}=\sqrt{1+x^2}$. 
Throughout this paper, the nonlinear term $N(u,\pa_x u)$ is always assumed to 
be a cubic homogeneous polynomial in  $(u,\cc{u}, \pa_x u, \cc{\pa_x u})$ with 
complex coefficients.  We will often write $u_x$ for $\pa_x u$.

From the perturbative point of view, 
cubic nonlinear Schr\"odinger equations in 
one space dimension are of special interest because the best possible 
decay in $L_x^2$ of general cubic nonlinear terms is $O(t^{-1})$, 
so the cubic nonlinearity 
must be regarded as a long-range perturbation. What we can expect 
for general cubic nonlinear Schr\"odinger equations in $\R$ 
is the lower estimate for the lifespan $T_{\eps}$ 
in the form $T_{\eps}\ge \exp(c/\eps^2)$ with some $c>0$, 
and this is best possible in general (see \cite{Kita} for an example of 
small data blow-up). More precise information on the lower bound 
is available under the restriction 
\begin{align}
 N(e^{i\theta},0) =e^{i \theta} N(1,0), \qquad \theta \in \R.
\label{weak_gi}
\end{align}
According to \cite{SagSu} (see also \cite{Su2}), 
if we assume \eqref{weak_gi} and the initial condition in \eqref{nls} 
is replaced by $u(0,x)=\eps \psi(x)$ with $\psi\in H^3\cap H^{2,1}$, 
then it holds that 
\begin{align}
 \liminf_{\eps \to +0} \eps^2 \log T_{\eps} 
\ge 
\frac{1}{\dis{2\sup_{\xi \in \R}(|\op{F}\psi(\xi)|^2 \imagpart \nu(\xi))}}
\label{lifespan}
\end{align}
with the convention $1/0=+\infty$,
where the function $\nu:\R\to \C$ is defined by 
\begin{align}\label{def_nu}
 \nu(\xi)=\frac{1}{2\pi i} \oint_{|z|=1} N(z,i\xi z) \frac{dz}{z^2}
\end{align}
and $\op{F}$ denotes the Fourier transform, i.e.,
\[
 \bigl(\op{F}\psi \bigr)(\xi)
 =
 \frac{1}{\sqrt{2\pi}} \int_{\R} e^{-iy\xi} \psi(y)\, dy
\]
for $\xi \in \R$.
Note that \eqref{weak_gi} excludes just the worst terms $u^3$, $|u|^2\cc{u}$, 
$\cc{u}^3$. As pointed out in \cite{HN3}, \cite{HN4}, \cite{HN5}, \cite{HN6}, 
\cite{HN7}, \cite{MP}, \cite{Nau}, etc., these three terms make the 
situation much more complicated. 
We do not intend to pursue this case here.  
We always assume \eqref{weak_gi} in what follows. 

In view of the right-hand side in \eqref{lifespan}, it may be natural 
to expect that the sign of $\imagpart \nu (\xi)$ has something to do with 
global behavior of small data solutions to \eqref{nls}. 
In fact, it has been pointed out in \cite{SagSu} that typical results on 
small data global existence and large-time asymptotic behavior for 
\eqref{nls} under \eqref{weak_gi}  can be summarized in terms of 
$\imagpart \nu(\xi)$ as follows. 
\begin{itemize}
\item[{\bf (i)}]
The small data global existence holds in $H^3\cap H^{2,1}$ 
under the condition 
\begin{align}
\imagpart \nu(\xi) \le 0,\quad \xi \in \R.
\tag{${\bf A}$}
\end{align}
\item[{\bf (ii)}]
If the inequality in (${\bf A}$) is replaced by the equality, i.e., 
\begin{align}
\imagpart \nu(\xi) = 0,\quad \xi \in \R, 
\tag{${\bf A}_0$}
\end{align}
then the  solution has 
a logarithmic phase correction in the asymptotic profile, 
i.e., it holds that 
\[
 u(t,x)= \frac{1}{\sqrt{t}} \alpha^+(x/t) 
 \exp\left( 
  \frac{ix^2}{2t} - i|\alpha^+(x/t)|^2 \realpart \nu (x/t) \log t
 \right)
 + 
o(t^{-1/2})
\]
as $t \to +\infty$ uniformly in $x \in \R$, 
where  $\alpha^+(\xi)$ is a suitable $\C$-valued function 
of $\xi \in \R$.
\item[{\bf (iii)}]
If the inequality in (${\bf A}$) is strict, i.e.,
\begin{align}
\sup_{\xi \in \R} \imagpart \nu(\xi)<0, 
\tag{${\bf A}_+$}
\end{align}
then the solution gains an additional logarithmic time decay
$\|u(t)\|_{L^{\infty}}  =O((t\log t)^{-1/2})$.
\end{itemize}
For more details on each case,  see the references cited in 
Section~1 of \cite{SagSu}. 
As for the large time behavior in the sense of $L_x^2$ under 
(${\bf A}$), 
it is not difficult to see that (${\bf A}_+$) 
implies $\dis{\lim_{t\to +\infty}\|u(t)\|_{L^2}= 0}$, whereas 
(${\bf A}_0$) implies 
$\dis{\lim_{t\to +\infty}\|u(t)\|_{L^2}\ne 0}$ for generic initial data 
of small amplitude. However, it is not clear whether $L^2$-decay occurs 
or not in the other cases
(even for a simple example such as $N(u,u_x)=-i|u_x|^2(u +u_x) +(u^3)_x$, 
for which we have
\[
\nu(\xi)
=
\frac{1}{2\pi i} \oint_{|z|=1} (-i\xi^2(1+i\xi) |z|^2z +3i\xi z^3)
\frac{dz}{z^2}
=
-i\xi^2 +\xi^3
\]
and $\imagpart \nu(\xi)=-\xi^2$). Despite the recent progress 
of studies on dissipative nonlinear Schr\"odinger equations 
(\cite{HLN}, \cite{HNS}, \cite{Hoshino}, \cite{JJL}, \cite{KLS}, 
\cite{Kim}, \cite{KitaLi}, \cite{KitaNak}, \cite{KitaShim}, 
\cite{LiS1}, \cite{LiS2}, \cite{Shim}, etc.), 
questions on decay/non-decay in $L_{x}^2$ without (${\bf A}_+$)
have not been addressed in the previous works 
except \cite{LNSS1} and \cite{LNSS2}.

The aim of this paper is to fill in the missing piece between 
(${\bf A}_+$) and (${\bf A}_0$), that is,  
to investigate $L^2$-decay property  of global solutions to \eqref{nls} 
under \eqref{weak_gi} and (${\bf A}$) without 
(${\bf A}_+$) and (${\bf A}_0$). Our main result is as follows.

\begin{thm}  \label{thm_main}
Suppose that $\eps =\|\varphi\|_{H^3\cap H^{2,1}}$ 
is sufficiently small. Assume that \eqref{weak_gi} and $({\bf A})$ are 
satisfied but $({\bf A}_0)$ is violated. 
Then, for any $\delta>0$, there exists a positive constant $C$ 
such that the global solution $u$ to \eqref{nls} satisfies
\[
\|u(t)\|_{L^2}\le \frac{C\eps}{(1+ \eps^2\log (t+2))^{1/4-\delta}}
\] 
for $t\ge 0$.
\end{thm}

\begin{rmk}
Under  \eqref{weak_gi} and (${\bf A}_+$), we can show 
the global solution to  \eqref{nls} has the stronger $L^2$-decay of 
order $O((\log t)^{-3/8+\delta})$ with arbitrarily small $\delta>0$ by 
the same method. 
For the detail, see Remark~\ref{rmk_final} below.
\end{rmk}

\begin{rmk} 
An analogous result for semilinear wave equations in $\R^2$ 
has been obtained in \cite{NST}, where the condition corresponding to 
(${\bf A}$) is called the Agemi condition.
\end{rmk}

\begin{rmk}\label{rmk_system}
In the case of systems, 
the situation is much more delicate than the single 
case. Detailed discussions on a weakly dissipative nonlinear Schr\"odinger 
system relevant to the present work can be found in \cite{LNSS1} and 
\cite{LNSS2} (see also \cite{NS} for a closely related work on a system of 
semilinear wave equations in $\R^2$).
\end{rmk}

\section{Proof}  \label{sec_proof}
The rest part of this paper is devoted to the proof of Theorem~\ref{thm_main}. 
The argument will be divided into four steps. 
\\

\noindent\underline{\bf Step 1:}\ 
We begin with the following elementary lemma, whose proof is skipped.
\begin{lem} \label{lem_key}
Let $p(\xi)$ be a real polynomial with $\deg p \le 3$. 
If $p(\xi) \ge 0$ for all $\xi \in \R$, then we have either of the following 
three assertions.
\begin{itemize}
\item[$({\bf a})$] 
$p(\xi)$ vanishes identically on $\R$.
\item[$({\bf b})$] 
$\dis{\inf_{\xi \in \R}p(\xi)>0}$.
\item[$({\bf c})$]
There exist $c_0>0$ and $\xi_0\in \R$ such that 
$p(\xi)=c_0(\xi-\xi_0)^2$.
\end{itemize}
\end{lem}

For $\nu(\xi)$ given by \eqref{def_nu}, we put $p(\xi)=-\imagpart\nu(\xi)$. 
Since we assume that (${\bf A}$) is satisfied but (${\bf A}_0$) is violated, 
we see that the case ({\bf a}) in Lemma~\ref{lem_key} is excluded. 
Note also that ({\bf b}) is equivalent to (${\bf A}_+$). 
Now, let us turn our attentions to the admissible range of the parameter 
$\theta$ for convergence of the integral 
\begin{align}
\label{integral}
 I_{\theta}=\int_{\R} \frac{d\xi}{p(\xi)^{\theta} \jb{\xi}^{4-4\theta}}
\end{align}
under ({\bf c}) or ({\bf b}). 
In the case ({\bf c}), we have 
\[
 I_{\theta}
=
 c_0^{-\theta} 
 \int_{\R} \frac{d\xi}{|\xi-\xi_0|^{2\theta} \jb{\xi}^{4-4\theta}}
<\infty
\]
for  $\theta<1/2$. 
In the case ({\bf b}), we have 
\[
 I_{\theta}
 \le
 \bigl(\inf_{\xi \in \R}p(\xi)\bigr)^{-\theta}
 \int_{\R} \frac{d\xi}{\jb{\xi}^{4-4\theta}}
<\infty
\]
for  $\theta<3/4$. \\

\noindent\underline{\bf Step 2:}\ 
Next we summarize the basic estimates for the global 
solution $u$ to \eqref{nls}. In what follows, we denote various positive 
constants by the same letter $C$ 
which may vary from one line to another.

First we write $\op{J}=x+it\pa_x$ and $\op{L}=i\pa_t+\frac{1}{2}\pa_x^2$. 
We note the important commutation relations
$[\pa_x, \op{J}]=1$, $[\op{L}, \op{J}]=0$. 
Next we set $\op{U}(t)=\exp(i\frac{t}{2}\pa_x^2)$ and 
$\alpha(t,\xi)=\op{F}\bigl[\op{U}(-t) u (t,\cdot)\bigr](\xi)$ 
for the solution $u$ to \eqref{nls}.
According to the previous works (\cite{HN2}, \cite{HNS}, \cite{SagSu}, etc.), 
we already know the following estimates. 
\begin{lem} \label{lem_apriori}
Let $\eps=\|\varphi\|_{H^{3}\cap H^{2,1}}$ be suitably small. 
Assume that \eqref{weak_gi} and $({\bf A})$ are fulfilled. Then, 
the solution $u$ to \eqref{nls} satisfies 
\begin{align}
 |\alpha(t,\xi)| \le \frac{C\eps}{\jb{\xi}^2}
\label{est_alpha}
\end{align}
for $t\ge 0$, $\xi \in \R$, and 
\begin{align}
\|u(t)\|_{H^3}+\|\op{J}u(t)\|_{H^2}\le C\eps (1+t)^{\gamma}
\label{est_sobolev}
\end{align}
for $t\ge 0$, where $0<\gamma<1/12$.
\end{lem}

The following lemma has been obtained in \cite{SagSu} (see also \cite{HN2}).
We write $\alpha_{\omega}(t,\xi)=\alpha(t,\xi/\omega)$ for 
$\omega\in \R\backslash\{0\}$.
\begin{lem} \label{lem_decomp}
Under the assumption \eqref{weak_gi}, we have 
\begin{align}
\op{F}  \op{U}(-t) (1-\pa_x^2) N(u,u_x) 
=&
(1+\xi^2)\frac{\nu(\xi)}{t} |\alpha|^2 \alpha 
+
\frac{\xi e^{i\frac{t}{3}\xi^2}}{t} \mu_{1}(\xi) \alpha_3^3 
\notag\\
&+
\frac{\xi  e^{i\frac{2t}{3}\xi^2} }{t} \mu_{2}(\xi) 
 \bigl(\cc{\alpha_{-3}}\bigr)^3
+
\frac{\xi e^{it\xi^2}}{t} \mu_{3}(\xi) |\alpha_{-1}|^2 \cc{\alpha_{-1}}  
+R, 
\label{crucial}
\end{align}
where $\nu(\xi)$ is given by \eqref{def_nu}, 
$\mu_{1}(\xi)$, $\mu_{2}(\xi)$, $\mu_{3}(\xi)$ are polynomials in 
$\xi$ of degree at most $4$, and $R(t,\xi)$ satisfies
\begin{align}
 |R(t,\xi)|
 \le
 \frac{ C}{t^{5/4}} \bigl(\|u(t)\|_{H^3}+\|\op{J}u(t)\|_{H^2} \bigr)^3
\label{est_R1}
\end{align}
for $t\ge 1$, $\xi \in \R$.
\end{lem}
For the proof, see Lemma~4.3 in \cite{SagSu}. 
By \eqref{est_sobolev} and \eqref{est_R1}, we have 
\begin{align}
 |R(t,\xi)|
 \le
 \frac{ C\eps^3}{t^{1+\kappa}}
\label{est_R2}
\end{align}
for $t\ge 1$, $\xi\in \R$, where $\kappa=1/4-3\gamma>0$.
This indicates that $R$ can be regarded as a remainder in \eqref{crucial}. 
We also observe that one $\xi$ pops up in front of the oscillating factors 
in \eqref{crucial}. This is the point where \eqref{weak_gi} plays a crucial 
role. 
As for the role of $\nu(\xi)$, the first term of the right-hand side in 
\eqref{crucial} tells us that $\nu(\xi)$ is responsible for the contribution 
from the gauge-invariant part in $N$.
\\

\noindent\underline{\bf Step 3:}\ 
We are going to make some reductions.
The goal in this step is to derive the ordinary differential equation 
\eqref{ode_profile} (with $\xi\in \R$ regarded as a parameter). 

Let $t\ge 2$ from now on. By the relation 
$\op{L}=\op{U}(t) i\pa_t \op{U}(-t)$ and Lemma~\ref{lem_decomp}, we have 
\begin{align}
 i\pa_t \alpha(t,\xi)
 &=
 \op{F} \op{U}(-t) \op{L}u
 \nonumber\\
 &=
 \jb{\xi}^{-2}\op{F} \op{U}(-t) (1-\pa_{x}^2)N(u,u_x)
 \nonumber\\
 &=
 \frac{\nu(\xi)}{t} |\alpha(t,\xi)|^2\alpha(t,\xi)
 + \eta(t,\xi)+ \jb{\xi}^{-2}R(t,\xi),
 \label{profile_eq}
\end{align}
where
\[
 \eta(t,\xi)
 =
 \frac{\xi e^{it\xi^2/3}}{t} 
 \frac{\mu_1(\xi)}{\jb{\xi}^2} \alpha_3^3
 + 
 \frac{\xi e^{i2t\xi^2/3}}{t} 
 \frac{\mu_2(\xi)}{\jb{\xi}^2} \cc{\alpha_{-3}}^3
 + 
 \frac{\xi e^{it\xi^2}}{t} 
 \frac{\mu_3(\xi)}{\jb{\xi}^2} |\alpha_{-1}|^2 \cc{\alpha_{-1}}.
\]
It follows from \eqref{est_alpha}, \eqref{est_R2} and \eqref{profile_eq} 
that
\begin{align*}
 |\pa_t \alpha (t,\xi)| 
 \le
  \frac{C\jb{\xi}^3}{t} \left(\frac{C\eps}{\jb{\xi}^2}\right)^3
 + 
 \frac{C\eps^3}{t^{1+\kappa}\jb{\xi}^2}
 \le 
 \frac{C\eps^3}{t \jb{\xi}^2}.
\end{align*}
Also, by using the identity
\begin{align*}
 \frac{\xi  e^{i\omega t\xi^2}}{t} f(t,\xi)
&= 
 \frac{\xi  \pa_t(te^{i\omega t\xi^2})}{t(1+i \omega t\xi^2)} f(t,\xi)\\
&=
 i \pa_t 
 \left(\frac{ -i\xi e^{i\omega t\xi^2}}{1+i \omega t\xi^2} f(t,\xi)\right) 
 -
 t e^{i \omega t\xi^2}
 \pa_t\left( \frac{ \xi f(t,\xi)}{t(1+i\omega t\xi^2)}  \right)
\end{align*}
and the inequality 
\[
 \sup_{\xi \in \R} 
 \frac{|\xi|^a }{|1+i \omega t\xi^2|}
 \le 
 \frac{C}{(|\omega|t)^{a/2}}
\]
for $0\le a \le 2$, we see that $\eta(t,\xi)$ can be split into 
\begin{align}
\eta=i\pa_t \sigma_1 +\sigma_2; 
\quad
 |\sigma_1(t,\xi)| \le \frac{C\eps^3}{t^{1/2}\jb{\xi}^4}, 
\quad
 |\sigma_2(t,\xi)| \le \frac{C\eps^3}{t^{3/2} \jb{\xi}^4}.
\label{decomp_eta}
\end{align}
With this $\sigma_1$, we set 
$\beta(t,\xi)=\alpha(t,\xi) -\sigma_1(t,\xi)$. 
Then it follows from \eqref{profile_eq} that  
\begin{align}\label{ode_profile}
 i\pa_t \beta(t,\xi) 
 =
 \frac{\nu(\xi)}{t} |\beta(t,\xi)|^2 \beta(t,\xi) +\rho(t,\xi),
\end{align}
where
\begin{align*}
 \rho(t,\xi)
 =&
 \frac{\nu(\xi)}{t}
 \biggl(|\alpha |^{2}\alpha    -|\beta|^2\beta \biggr)
  +
 \sigma_2
 + 
 \jb{\xi}^{-2}R\\
=&
 \frac{\nu(\xi)}{t}
 \biggl( 
   2|\alpha|^{2}\sigma_{1}
   +
    \alpha^{2}\overline{\sigma_{1}}
   -
   2\alpha|\sigma_{1}|^{2}
   - 
   \overline{\alpha}\sigma_{1}^{2}
   +
   |\sigma_{1}|^{2}\sigma_{1}
 \biggr)
  +
 \sigma_2
 + 
 \jb{\xi}^{-2}R.
\end{align*}
By \eqref{est_alpha}, \eqref{est_R2} and \eqref{decomp_eta}, we have 
\begin{align*}
 |\rho(t,\xi)|
 &\le
  \frac{C\jb{\xi}^3}{t} 
 \left(
   \frac{C\eps}{\jb{\xi}^2}  +  \frac{C\eps^3}{t^{1/2} \jb{\xi}^4}
 \right)^2
 \frac{C\eps^3}{t^{1/2}\jb{\xi}^4}
 + 
 \frac{C\eps^3}{t^{3/2} \jb{\xi}^4}
 + 
 \frac{C\eps^3}{t^{1+\kappa}\jb{\xi}^2}\\
 &\le
 \frac{C \eps^3}{t^{1+\kappa}\jb{\xi}^2}. 
\end{align*}
Remember that $0<\kappa<1/4$.

Roughly speaking, what we have seen so far is that the solution 
$u$ to \eqref{nls} under \eqref{weak_gi} can be expressed as 
\[
 u=\op{U}(t)\op{F}^{-1}\beta+\cdots
\]
with 
\[ i\pa_t \beta =\frac{\nu(\xi)}{t}|\beta|^2 \beta+\cdots,
\]
where the terms ``$+\cdots$" are expected to be  harmless. 
By this reason it would be fair to call \eqref{ode_profile} 
the {\em profile equation} associated with \eqref{nls} under 
\eqref{weak_gi}. The original idea of this reduction is 
due to Hayashi--Naumkin \cite{HN1}.\\

\noindent\underline{\bf Final step:}\ 
Before going further, let us recall the following useful lemma 
due to Matsumura. 
\begin{lem}\label{lem_Matsumura}
Let $C_0>0$, $C_1\ge 0$, $q>1$ and $s>1$. 
Suppose that a function $\Phi(t)$ satisfies 
\begin{align*}
\frac{d\Phi}{dt}(t)
\le -\frac{C_0}{t}\left|\Phi(t)\right|^q + \frac{C_1}{t^{s}} 
\end{align*}
for $t\ge 2$. Then we have
\begin{align*}
\Phi(t)\le \frac{C_2}{(\log t)^{q^{*}-1}}
\end{align*}
for $t\ge 2$, where $q^{*}$ is the H\"older conjugate of $q$ (i.e., 
$1/q+1/q^{*}=1$), and
\begin{align*}
C_2=\frac{1}{\log 2}\left( (\log 2)^{q^{*}}\Phi(2) 
  + C_1\int_{2}^{\infty}\frac{(\log \tau)^{q^{*}}}{\tau^{s}}d\tau \right) 
  + \left( \frac{q^{*}}{q C_0} \right)^{q^{*}-1}.
\end{align*}
\end{lem}
For the proof, see Lemma~4.1 in \cite{KMatsS}. 
At last, we are in a position to reach the conclusion. 
We set $\Phi(t,\xi)=p(\xi) |\beta(t,\xi)|^2$ 
with $p(\xi)=-\imagpart \nu(\xi)$.
Note that $\Phi(t,\xi) \ge 0$ by (\textbf{A}). 
It follows from \eqref{ode_profile} that
\begin{align*}
 \pa_t\Phi(t,\xi) 
 &=
 2p(\xi)\imagpart\bigl(\cc{\beta(t,\xi)}i\pa_t\beta(t,\xi)\bigr)
 \\
 &=
 2p(\xi) 
 \left( 
  \frac{\imagpart\nu(\xi)}{t} |\beta(t,\xi)|^4 
  +
  \imagpart\bigl(\cc{\beta(t,\xi)}\rho(t,\xi) \bigr) 
 \right)
 \\
 &\le
  -\frac{2p(\xi)^2}{t}|\beta(t,\xi)|^4
 + 
 C\jb{\xi}^3 \frac{C\eps}{\jb{\xi}^2} \frac{C\eps^3}{t^{1+\kappa}\jb{\xi}^2}
 \\
 &\le 
 -\frac{2}{t}\Phi(t,\xi)^2+\frac{C\eps^4}{t^{1+\kappa}\jb{\xi}},
\end{align*}
where $\kappa\in (0,1/4)$.
We also note that \eqref{est_alpha} yields
\[
 \Phi(2,\xi) 
 \le 
  C\jb{\xi}^3 \left(\frac{C\eps}{\jb{\xi}^2}\right)^2
 \le
  \frac{C\eps^2}{\jb{\xi}}.
\]
Therefore we can apply Lemma~\ref{lem_Matsumura} with $q=2$ and $s=1+\kappa$ 
to obtain 
\[
 0\le \Phi(t,\xi) \le \frac{C}{\log t},
\]
whence
\begin{align*}
 |\alpha(t,\xi)| 
&\le 
\sqrt{\frac{\Phi(t,\xi)}{p(\xi)}} +|\sigma_1(t,\xi)|
\\
&\le 
\frac{C}{\sqrt{p(\xi)\log t}} 
\left(
 1 + \eps^3 \frac{\sqrt{p(\xi)}}{\jb{\xi}^4} \sqrt{\frac{\log t}{t}}
\right)
\\
&\le 
\frac{C\eps}{\sqrt{p(\xi)\eps^2 \log t}}.
\end{align*}
Interpolating this with \eqref{est_alpha}, we deduce that 
\begin{align}
 |\alpha(t,\xi)| 
\le 
\frac {C\eps}{(\eps^2 \log t)^{\theta/2}}
\frac {1}{p(\xi)^{\theta/2}\jb{\xi}^{2-2\theta}}
\label{decay_main}
\end{align}
for $\theta \in [0,1]$. By the $L^2$-unitarity of $\op{U}(t)$ and $\op{F}$, 
we have 
\begin{align}
\|u(t)\|_{L^2}^2
=
\|\alpha(t)\|_{L^2}^2
\le
\frac
 {C\eps^2}{(\eps^2\log t)^{\theta}}I_{\theta}
\label{L2_est}
\end{align}
for $0\leq \theta <\frac{1}{2}$, 
where $I_{\theta}$ is given by \eqref{integral}. Therefore 
we can take $\theta=1/2-2\delta$ with $\delta>0$ to see that 
\[
\|u(t)\|_{L^2}\le \frac
 {C\eps}{(\eps^2\log t)^{1/4-\delta}}.
\]
Also we obtain $\|u(t)\|_{L^2} \le C\eps$ by taking $\theta =0$ 
in \eqref{L2_est}. Piecing them together, we arrive at the desired estimate. 
\qed\\

\begin{rmk}\label{rmk_final}
Under \eqref{weak_gi} and the stronger condition (${\bf A}_+$), it is possible 
to choose 
$\theta=3/4-2\delta$ in \eqref{L2_est} because (${\bf A}_+$) implies 
({\bf b}) in Lemma~\ref{lem_key} 
and thus the admissible range for $\theta$ in \eqref{L2_est} becomes 
$0\leq \theta <\frac{3}{4}$. 
That is the reason why $\|u(t)\|_{L^2}$ decays like 
$O((\log t)^{-3/8+\delta})$ under (${\bf A}_+$). 
It is not certain whether this rate is the best or not. 
Indeed, it is possible to improve the exponent from 
$-3/8+\delta$ to $-1/2$ if there exists a positive constant $C_*$ 
such that 
\begin{align}
\imagpart \nu(\xi) \le -C_*\jb{\xi}^2, \quad \xi \in \R 
\tag{${\bf A}_{++}$}
\end{align}
(cf. Theorem~2.3 in \cite{LiS1}). A typical example of $N$ satisfying 
(${\bf A}_{++}$) is $-i|u+u_x|^2u$.

It may be an interesting problem to specify the optimal $L^2$-decay rates
for the solutions to \eqref{nls} under \eqref{weak_gi} and (${\bf A}$) 
(with or without (${\bf A}_+$)).
\end{rmk}

\medskip
\subsection*{Acknowledgments}
The work of H.~S. is supported by Grant-in-Aid for Scientific Research (C) 
(No.~17K05322), JSPS. 


\end{document}